\documentclass{amsart}
\input amssym.def\input amssym

\newtheorem{thm}{Theorem}
 \newtheorem{cor}[thm]{Corollary}
 \newtheorem{lem}[thm]{Lemma}
 
 \newtheorem{defn}[thm]{Definition}

\newcommand{\ma}{\mathcal{A}}

\newcommand{\mm}{\mathcal{M}}

\begin{document}

\title{Additive preserving rank one maps on Hilbert $C^\ast$-modules }

\author{Bin Meng}
\address{College of Science, Nanjing University of Aeronautics and Astronautics, Nanjing 210016,  People's Republic
 of China}

 \email{b.meng@nuaa.edu.cn}


\date{}


\begin{abstract}In this paper, we characterize a class of additive
maps on Hilbert $C^\ast$-modules which maps a "rank one" adjointable
operators to another rank one operators.
\end{abstract}

\maketitle

\section{Introduction and preliminaries}
The problem of determining linear maps $\Phi$ on \textbf{B}(X)
preserving certain properties has attracted attention of many
mathematicians in resent decades. They have been devoted to the
study of linear maps preserving spectrum, rank, nilpotency, etc.

Rank preserving problem is a basic problem in the study of linear
preserver problem. Rank preserving linear maps have been studied
intensively by Hou in \cite{hou1}. In \cite{kuz,oml} the authors
used very elegant arguments to completely describe additive mappings
preserving (or decreasing) rank one.

In our study of free probability theory, we find that module maps on
Hilbert $C^\ast-$module preserving certain properties are also
important \cite{meng1,meng2}, and thus the study of the modules
preserving problem becomes attractive. In \cite{meng3}, a complete
description of modules maps which preserving rank one modular
operators is given. Naturally in the present paper we consider
additive maps on modules preserving rank one, which much more
complicated than the module maps case.

A Hilbert $C^\ast-$module over a $C^\ast-$module over a
$C^\ast-$algebra $\ma$ is a left $\ma-$module $\mm$ equipped with an
$\ma-$valued inner product $\langle, \rangle$ which is $\ma-$linear
in the first and conjugate $\ma-$linear in the second variable such
that $\mm$ is a Banach space with the norm $\|v\|_M=\|\langle
v,v\rangle\|^{\frac{1}{2}}$, $\forall v\in\mm$. Hilbert
$C^\ast-$modules first appeared in the work of Kaplansky \cite{kap},
who used them to prove that derivations of type I $AW^\ast-$algebras
are inner. Now a good text book about Hilbert $C^\ast-$module is
\cite{lan}. In this paper we mainly consider Hilbert $\ma-$module
$H\otimes\ma$ (or denoted by $H_\ma$), where $H$ is a separable
infinite dimentional Hilbert space and $\ma$ is a unital
$C^\ast-$algebra. $H_\ma$ play a special role in the theory of
Hilbert $C^\ast-$modules (see \cite{lan}). Obviously, $H_\ma$ is
countably generated and possesses orthonormal basis $\{e_i\otimes
1\}$, which $\{e_i\}$ is an orthonormal basis in $H$.

We introduce a class of module maps which is analogues to rank one
operators on a Hilbert space. For all $x,y\in H_\ma$, define
$\theta_{x,y}:H_\ma\rightarrow H_\ma$ by $\theta_{x,y}(\xi):=\langle
\xi,y\rangle x$, for every $\xi\in H_\ma$. Note that $\theta_{x,y}$
is quite different from rank one linear operators on Hilbert space.
For instance, we cannot infer $x=0$ or $y=0$ from $\theta_{x,y}=0$.
We denote $\{\sum\limits_{i=1}^n \alpha_i \theta_{x_i,y_i},\forall
n\in \Bbb{N}\}$ by $\mathcal{F}(H_\ma)$.

In the present paper, we will describe the additive map
$\Phi:\mathcal{F}(H_\ma)\rightarrow \mathcal{F}(H_\ma)$ which maps
$\theta_{x,y}$ to some $\theta_{s,t}$. The methods are analogues to
those in \cite{hou1,kuz,oml} but much more complicated.

\section{Definitions and lemmas}
In this section we mainly introduce some definitions and prove some
lemmas.

\begin{defn}
\cite{meng3} Let $\varepsilon$ be a orthonormal basis in $H_\ma$.
$x\neq 0$ in $H_\ma$ will be called coordinatly invertible if for
all $e\in\varepsilon$, $\langle e,x\rangle$ is invertible in $\ma$
unless $\langle e,x\rangle=0$.
\end{defn}

We denote the set of all the coordinatly invertible elements in
$H_\ma$ by $CI(H_\ma$ or CI for short.

Coordinatly invertible elements in Hilbert $C^\ast-$module are
analogues to elements in linear space to some extents.
\begin{lem}
\cite{meng3} Supposing $y\in CI$ and $\theta_{x,y}=0$ then $x=0$.
\end{lem}
\begin{lem}
\cite{meng3} Let $\mm$ be Hilbert $\ma-$modules, where $\ma$ is a
unital $C^\ast-$algebra, and let $\phi,\sigma:\mm\rightarrow\ma$ be
$\ma-$linear operators. Suppose that $\sigma$ vanishes on the kenel
of $\phi$. Then there exists a $b\in\ma$ such that $\sigma=\phi\cdot
b$.
\end{lem}

\begin{cor}
\cite{meng3}Let $g_1,g_2\in\mm$. If for all $x\in\mm$, $\langle x,
g_1\rangle=0$ implying $\langle x,g_2\rangle=0$, then there is a
$a\in\ma$ such that $g_2=ag_1$.
\end{cor}

The following Lemma 5, which will be used frequently, has been
obtained in \cite{meng3}. Nevertheless we give its proof for the
sake of completeness.

\begin{lem}
Let $\ma$ be a unital $C^\ast-$algebra, $x_1,x_2\in\mm$, $g_1,g_2\in
CI$ satisfying $\theta_{x_1,g_1}+\theta_{x_2,g_2}=\theta_{x_3,g_3}$.
Then at least one of the following is true:\\
(i)there exists a invertible $\alpha_1\in\ma$ such that $g_1=\alpha_1g_2$;\\
(ii) there are $\beta_1,\beta_2\in\ma$ such that $x_1=\beta_1 x_3$,
$x_2=\beta_2x_3$.
\end{lem}

\begin{proof}
We will complete the proof by considering the following four cases.

Case 1. For all $\xi\in H_\ma$, $\langle \xi, g_2\rangle=0$ implying
$\langle\xi, g_1\rangle=0$. From Corollary 4, there exists
$\alpha_1\in\ma$ such that $g_1=\alpha_1 g_2$. Furthermore since
$g-1,g_2\in CI$, we infer that $\alpha_1\in\ma$ is invertible.

Case 2. For all $\xi\in H_\ma$, $\langle\xi,g_1\rangle=0$ implying
$\langle\xi, g_2\rangle=0$. Still from Corollary 4, there is a
$\alpha_2\in\ma$ such that $g_2=\alpha_2 g_1$ and $\alpha_2$ is
invertible.

Case 3. There exists a $\xi_0\in H_\ma$ such that
$\langle\xi_0,g_2\rangle=0$ but $\langle \xi_0, g_1\rangle\neq 0$.
We can find $e\in\varepsilon$ such that $\langle e,g_2\rangle=0 $
but $\langle e,g_1\rangle\neq 0$. Then from $\langle e,g_1\rangle
x_1+\langle e,g_2\rangle x_2=\langle e,g_3\rangle x_3$, it follows
$\langle e,g_1\rangle x_1=\langle e,g_3\rangle x_3$. Since $g_1\in
CI$, we have $x_1=\langle e, g_1\rangle ^{-1}\langle e,g_3\rangle
x_3$. We put $\beta_1=\langle e, g_1\rangle^{-1}\langle e,
g_3\rangle$ and get $\theta_{\beta
x_3,g_1}+\theta_{x_2,g_2}=\theta_{x_3,g_3}$. Thus
$\theta_{x_2,g_2}=\theta_{x_3,g_3-\beta_1^\ast g_1}$. Now choosing a
$e'\in\varepsilon$, we get $\langle e',g_2\rangle x_2=\langle
e',g_3-\beta_1^\ast g_1\rangle x_3$ and thus $x_2=\langle
e',g_2\rangle^{-1}\langle e',g_3-\beta_1^\ast g_1\rangle x_3$.
Putting $\beta_2=\langle e',g_2\rangle ^{-1}\langle
e',g_3-\beta_1^\ast g_1\rangle$ then we obtain (ii).

Case 4. There exists $\xi_0\in H_\ma$ such that $\langle \xi_0,
g_1\rangle=0$ but $\langle \xi_0,g_2\rangle\neq 0$. Similar to Case
3, we get (ii) again.
\end{proof}

\begin{cor}
With the notations in the above lemma, suppose $g_1\neq \alpha g_2$,
for all $\alpha\in\ma$ and $g_3\in CI$. Then there exist
$\beta_1,\beta_2\in \ma$ which are invertible such that $x_1=\beta_1
x_3$, $x_2=\beta_2 x_3$.
\end{cor}
\begin{proof}
Denote $\{x|\langle x, g_i\rangle\}$ by $ker g_i$, $i=1,2$. Since
$g_1\neq \alpha g_2$, $g_2\neq \beta g_1$, we have $ker
g_1\nsubseteq ker g_2$, $ker g_2\nsubseteq ker g_1$. So there exists
$e_1\in H_\ma$, such that $\langle _1,g_1\rangle\neq 0$ but $\langle
e_1,g_2\rangle=0$. Then $\langle e_1,g_1\rangle x_1+\langle
e_1,g_2\rangle x_2=\langle e_1,g_3\rangle x_3$, i.e. $\langle
e_1,g_1\rangle x_1=\langle e_1,g_3\rangle x_3$ and thus $x_1=\langle
e_1,g_1\rangle^{-1}\langle e_1,g_3\rangle x_3$. Putting $\langle
e_1,g_1\rangle ^{-1}\langle e,g_3\rangle=\beta_1$ is invertible.

Similarly there exits a $e_2\in H_\ma$ such that $\langle
e_2,g_1\rangle=0$ but $\langle e_2,g-2\rangle\neq 0$. Then
$x_2=\langle e_2,g_2\rangle^{-1}\langle e_2,g_3\rangle x_3$. Putting
$\beta_2=\langle e_2,g_2\rangle^{-1}\langle e_2,g_3\rangle$.
\end{proof}

\begin{defn}
$\Phi:\mathcal{F}(H_\ma)\rightarrow \mathcal(H_\ma)$ is a map. If
for any $x\in H_\ma$, $y\in CI$, there are $s\in H_\ma$, $t\in CI$
such that $\Phi(\theta_{x,y})=\theta_{s,t}$. Then we call $\Phi$ is
rank one decreasing. If $x\neq 0$ implying $s\neq 0$ then $\Phi$
will be called rank one preserving.
\end{defn}
\begin{defn}
$\Phi:\mathcal{F}(H_\ma)\rightarrow \mathcal{F}(H_\ma)$ is an
additive map and for arbitrary $\theta_{x,y}$,
$\Phi(\lambda\theta_{x,y})=\tau_{x,y}(\lambda)\Phi(\theta_{x,y})$
where $\tau_{x,y}:\ma\rightarrow\ma$ is a map. Then $\Phi$ will be
called a locally quasi-modular map.

If there exists a $\tau:\ma\rightarrow\ma $ is a map such that
$\Phi(\lambda T)=\tau(\lambda)\Phi(T)$, for all $T\in
\mathcal{F}(H_\ma)$ then $\Phi$ will be called a $\tau-$
quasi-modular map.
\end{defn}

\begin{lem}
Let $\ma$ be a unital commutative $C^\ast-$algebra , let $A,B$ be
injective $\tau-$quasi-modular continuous maps on $H_\ma$ with
$dim_\ma(A(H_\ma))\geq 2$ and suppose $\tau$ is surjective. There is
$\lambda_x\in\ma$, such that $Bx=\lambda_x Ax$. Then $B=\lambda A$
for some $\lambda\in\ma$.
\end{lem}
\begin{proof}
There are $x_1\in CI$, $x_2\in H_\ma$ such that $\langle
Ax_1,Ax_2\rangle=0$. From the assumption there exist
$\lambda_1,\lambda_2,\lambda_3\in\ma$ such that
$Bx_1=\lambda_1Ax_1$, $Bx_2=\lambda_2 Ax_2$ and
$B(x_1+x_2)=\lambda_3 A(x_1+x_2)$. Therefore,
$(\lambda_1-\lambda_3)Ax_1+(\lambda_2-\lambda_3)Ax_2=0$ and
\begin{eqnarray*}
&&\langle (\lambda_1-\lambda_3)Ax_1+(\lambda_2-\lambda_3)Ax_2,
(\lambda_1-\lambda_3)Ax_1+(\lambda_2-\lambda_3)Ax_2\rangle\\
&=&(\lambda_1-\lambda_3)\langle Ax_1,Ax_1\rangle
(\lambda_1-\lambda_3)^\ast+(\lambda_2-\lambda_3)\langle
Ax_2,Ax_2\rangle (\lambda_2-\lambda_3)^\ast\\
&=& 0.
\end{eqnarray*}
It follows that $(\lambda_1-\lambda_3)\langle Ax_1,Ax_1\rangle
(\lambda_1-\lambda_3)^\ast=0$ and $(\lambda_2-\lambda_3)\langle
Ax_2,Ax_2\rangle (\lambda_2-\lambda_3)^\ast=0$. Furthermore
$(\lambda_1-\lambda_3)Ax_1=(\lambda_2-\lambda_3)Ax_2=0$. We infer
$Bx_1=\lambda_1 Ax_1=\lambda_3Ax_1$ and $Bx_2=\lambda_2
Ax_2=\lambda_3 Ax_2$.

For every $x\in H_\ma$ such that $Ax,Ax_1,Ax_2$ is an orthogonal
set, we claim that there is a $\lambda\in\ma $ such that $Bx=\lambda
Ax$. In fact, $Bx=\lambda_x Ax=\lambda_{x+x_1}Ax$ and
$Bx_1=\lambda_1Ax_1=\lambda_{x+x_1}Ax_1=\lambda_3 Ax_1$. It follows
that $\lambda_{x+x_1}x_1=\lambda_3 x_1$ from $A$ is injective. And
then $\lambda_{x+x_1}=\lambda_3$. Thus we obtain $Bx=\lambda_3 Ax$.

$A(H_\ma)$ is submodule of $H_\ma$ since $\tau$ is surjective. Let
$Ax_1,Ax_2,\cdots$ be a orthogonal basis in $A(H_\ma)$. For
arbitrary $x\in H_\ma$, $Ax=\alpha_1 Ax_1+\alpha_2Ax_2+\cdots$.
Since $\tau$ is surjective, there exist
$\lambda_1,\lambda_2,\cdots,$ such that
$\alpha_i=\tau(\lambda_i),i=1,2,\cdots$. Therefore
$Ax=A(\lambda_1x_1+\lambda_2x_2+\cdots)$ and
$x=\lambda_1x_1+\lambda_2x_2+\cdots$. We infer there exists a
$\lambda\in\ma$ such that
\begin{eqnarray*}
Bx&=&B(\lambda_1x_1+\lambda_2x_2+\cdots)\\
&=&\alpha_1 B(x_1)+\alpha_2B(x_2)+\cdots\\
&=&\alpha_1\lambda Ax_1+\alpha_2\lambda Ax_2+\cdots\\
&=&\lambda A(\lambda_1x_1+\lambda_2x_2+\cdots)\\
&=&\lambda Ax
\end{eqnarray*}
\end{proof}

Now we introduce some notations. $L_x:=\{\theta_{x,g}\mid g\in
H_\ma\}$, $R_f:=\{\theta_{x,f}\mid x\in H_\ma\}$,
$L_x^{CI}:=\{\theta_{x,g}\mid g\in CI\}$,
$R_f^{CI}:=\{\theta_{y,f}\mid y\in CI\}$.

The following lemma play important roles in this paper.
\begin{lem}
Let $\Phi$ be a additive rank one preserving map.  For every $x\in
H_\ma$ there exists either a $y\in H_\ma$ such that
$\Phi(L_x^{CI})\subseteq L_y^{CI}$ or an $f\in CI$ such that
$\Phi(L_x^{CI})\subseteq R_f$.
\end{lem}
\begin{proof}
Assume that there exists a $x_0\in H_\ma$ such that
$\Phi(L_{x_0}^{CI}\subseteq L_x^{CI}$, $\Phi(L_{x_0}^{CI})\subseteq
R_f$, for all $x\in H_\ma, f\in CI$. Then there are $f_1,f_2\in CI$
such that $\Phi(\theta_{x_0,f_1})=\theta_{x_1,g_1}$,
$\Phi(\theta_{x_0,f_2})=\theta_{x_1,g_2}$, where $x_1\neq \alpha_1
x$, $x_2\neq \alpha_2 x$, for all $x\in
H_\ma,\alpha_1,\alpha_2\in\ma$ and $g_1\neq \beta_1 g_2$,
$g_2\neq\beta_2 g_1$, for all $\alpha,\beta\in\ma$. Since $\Phi$ is
rank one preserving, $\Phi(\theta_{x_0,f_1+f_2})=\theta_{x_3,g_3}$
for some $x_3\in H_\ma,g_3\in CI$. On the other hand, from Lemma 5,
we get $\theta_{x_1,g_1}+\theta_{x_2,g_2}\neq \theta_{x_3,g_3}$.
Thus we reach a contradiction.
\end{proof}
\begin{lem}
At least one of the following is true

(1)for all $x\in H_\ma$, there exits $y\in H_\ma$ such that
$\Phi(L_x^{CI})\subseteq L_y^{CI}$;

(2)for all $x\in H_\ma$, there exits $f\in CI$ such that
$\Phi(L_x^{CI}\subseteq R_f$
\end{lem}

\begin{proof}
If there are $Dim_\ma \Phi(L_x^{CI})=1$ ($\forall f_1,f_2\in CI$,
$\Phi (\theta_{x,f_1})=\theta_{y,g_1}$,
$\Phi(\theta_{x,f_2})=\theta_{y,g_2}$ implying $g_1=\alpha g_2$.)
then both (1) and (2) hold.

Now we consider $Dim_\ma \Phi(L_x^{CI})\geq 2,Dim_\ma \Phi(R_f)\geq
2$. Otherwise suppose, to reach a contradicion, that
$\Phi(L_{x_0}^{CI})\subseteq L_{y_0}^{CI}$,
$\Phi(L_{x_1}^{CI})\subseteq R_{g_1}$, where $x_0,x_1\in H_\ma$,
$y_0\in H_\ma$, $g_1\in CI$. Since $Dim_\ma \Phi(L_x^{CI})\geq 2$,
there exists a $h\in CI$, such that
$\Phi(\theta_{x_0,h})=\theta_{y_0,g}$, $g\in CI$ with $g\neq \alpha
g_1$, $\forall \alpha\in\ma$. Since $Dim_\ma \Phi(R_f)\geq 2$, we
know there exist $k\in H_\ma$, such that
$\Phi(\theta_{x_1,k})=\theta_{z,g_1}$ with $z\neq \beta x$, $y_0\neq
\beta x$, $\forall x\in H_\ma$, $\beta\in\ma$. Letting $m\in CI$
such that $\Phi(\theta_{x_0,k})=\theta_{y_0,m}$, then
$\Phi(\theta_{x_0+x_1,k})=\theta_{z,g_1}+\theta_{y_0,m}$. Since
$z\neq \alpha x,y_0\neq \beta x$, $\forall x\in
H_\ma,\alpha,\beta\in\ma$, there exists a $\lambda\in\ma$ such that
$m=\lambda g_1$ and
$\Phi(\theta_{x_0,k})=\theta_{y_0,m}=\theta_{y_0,\lambda g_1}$. On
the other hand, $\Phi(\theta_{x_0,h})=\theta_{y_0,g}$,
$\Phi(\theta_{x_1,h})=\theta_{l,g_1}$, where $h,g,g_1\in CI$. From
$g_1\neq \alpha g$, for all $\alpha\in\ma$ and
$\Phi(\theta_{x_0+x_1,h})=\theta_{y_0,g}+\theta_{l,g_1}$, we infer
that there are $\alpha_0,\beta_0$ which are invertible such that
$y_0=\alpha_0 y_0'$, $l=\beta_0 y_0'$. Therefore,
$\Phi(\theta_{x_1,g})=\theta_{\beta_0 y_0',g_1}$,
$\Phi(\theta_{x_0,h})=\theta_{\alpha_1 y_0',g}$ and
\begin{eqnarray*}
&&\Phi(\theta_{x_0+x_1,h+k})=\Phi(\theta_{x_0,h})+\Phi(\theta_{x_0,k})+\Phi(\theta_{x_1,h})
+\Phi(\theta_{x_1,k})\\
&=&\theta_{\alpha_0y_0',g}+\theta_{\alpha_0y_0',\lambda
g_1}+\theta_{\beta_0 y_0',g_1}+\theta_{z,g_1}\\
&=&\theta_{\alpha_0y_0',g}+\theta_{(\lambda^\ast
\alpha_0y_0'+\beta_0y_0'+z),g_1}
\end{eqnarray*}
Since $g\neq \alpha g_1,z\neq\beta x, y_0\neq \gamma x$, for all
$\alpha,\beta,\gamma\in\ma$, we know $\Phi(\theta_{x_0+x_1,h+k})$ is
not rank one which contradicts to the assumption of $\Phi$.
\end{proof}

\begin{cor}
 Then one of the
following is true

(i)for all $f\in H_\ma$ there exists a $f\in H_\ma $ with
$\Phi(R_f^{CI})\subseteq R_g^{CI}$;

(ii)for all $f\in H_\ma$ there exists $y\in CI$ with
$\Phi(R_f^{CI})\subseteq L_y$
\end{cor}
\begin{lem}
Suppose $Im\Phi$ is neither contained in any $L_y$ or contained in
any $R_g$ and $Dim_\ma\Phi(L_x^{CI})\geq 2$. Then

(i)If for all $x\in H_\ma$ there exists a $y\in H_\ma$ such that
$\Phi(L_x^{CI}$, then $\Phi(R_f^{CI})\subseteq R_g^{CI}$;

(ii)If for all $x\in H_\ma$ there exists a $y\in CI$ such that
$\Phi(L_x^{CI})\subseteq R_y$. Then for all $f\in H_\ma$ there
exists a $z\in CI$ such that $\Phi(R_f^{CI})\subseteq L_z$.
\end{lem}
\begin{proof}
We only prove (i). Assume, to reach a contradiction, that we have
simulatously  $\Phi(L_x^{CI})\subseteq L_y^{CI}$ and
$\Phi(R_f^{CI})\subseteq L_z$, for some $z\in CI$. Since $Dim_\ma
\Phi(L_x^{CI})\geq 2$ there are $k',k'',g',g''\in CI$ with $g'\neq
\alpha g''$ such that $\Phi(\theta_{x,k'})=\theta_{y,g'}$,
$\Phi(\theta_{x,k''})=\theta_{y,g''}$.

Since $Im\Phi$ is not contained in any $L_g$, there are $x_1\in
H_\ma$, $k_1\in CI$ with $\Phi(\theta_{x_1,k_1})=\theta_{y_1,g_1}$
such that $y_1\neq \alpha h$, $y\neq \beta h$, for all
$\alpha,\beta\in \ma, h\in H_\ma$. Consequently
$\Phi(L_{x_1}^{CI})\subseteq L_{y_1}^{CI}$,
$\Phi(\theta_{x_1,k'})=\theta_{y_1,u}$. Therefore,
$\Phi(\theta_{x+x_1,k'})=\theta_{y,g'}+\theta_{y_1,u}$ and then
$u=\lambda g'$ for some invertible $\lambda\in\ma$ i.e.
$\Phi(\theta_{x_1,k'})=\theta_{\lambda^\ast y_1,g'}$ and
$\Phi(R_{k'}^{CI})\subseteq R_{g'}$. It follows
$\Phi(R_{k'}^{CI})\nsubseteq L_z$ which implying $\Phi(R_f)\subseteq
R_g$. So we reach a contradiction.
\end{proof}

\begin{lem}
Suppose $Im\Phi$ is not contained in any $L_y$ nor $R_f$ and
$Dim_\ma\Phi(L_x^{CI})\geq 2$. If for all $x\in CI$, there exists a
$y'\in H_\ma$ such that $\Phi(L_x^{CI})\subseteq L_{y'}^{CI}$, then
there exists a $y\in CI$ such that $\Phi(L_x^{CI})\subseteq
L_y^{CI}$.
\end{lem}

\begin{proof}
It follows from $\Phi(L_x^{CI})\subseteq L_{y'}^{CI}$ that
$\Phi(L_x^{CI})\nsubseteq R_f$, for all $f\in CI$. If the conclusion
of the Lemma were wrong, we will have $\Phi(L_x^{CI})\nsubseteq
L_y^{CI}$, for all $y\in CI$. Thus there exist $g_1\neq \alpha
g_2,y_1\neq \alpha y,y_2\neq\beta y,y_1,y_2\in CI$, for all $y\in
CI$ with $\Phi(\theta_{x,f_1})=\theta_{y_1,g_1}$,
$\Phi(\theta_{x,f_2})=\theta_{y_2,g_2}$. We claim that
$y_1\neq\alpha x,y_2\neq\beta x$, for all $x\in H_\ma$. If not,
there exists a $y_0'\in H_\ma$ such that $y_1\alpha_0
y_0',y_2=\beta_0y_0'$. Then it follows that $\alpha_0,\beta_0$ are
invertible and $y_0'$ is coordinately invertible which contradicting
to $\Phi(L_x^{CI})\nsubseteq L_y^{CI}$, for all $y\in CI$. So
$y_1\neq \alpha x,y_2\neq \beta x$ but this this contradicting to
that $\Phi(\theta_{x,f_1+f_2})$ is rank one.
\end{proof}

\begin{lem}
If $\Phi(L_x^{CI})\subseteq L_y^{CI}$, then $\Phi(L_x)\subseteq
L_y$.
\end{lem}
\begin{proof}
For all $f\in H_\ma$, $f=\sum\limits_i \alpha_i f_i$, $f_i\in CI$,
we have $\Phi(\theta_{x,f})=\Phi(\theta_{x,\sum\limits_i
\alpha_if_i})=\sum\limits_i \alpha_i^\ast
\Phi(\theta_{x,f_i})=\sum\limits_i\alpha_i^\ast\theta_{y,g_i}=\theta_{y,\sum\limits_ig_i}$.
Then $\Phi(L_x)\subseteq L_y$.
\end{proof}
Similarly we can prove
\begin{lem}
If $\Phi(L_x^{CI})\subseteq R_f$, for some $f\in CI$, then
$\Phi(L_x)\subseteq R_f$.
\end{lem}

\begin{lem}
For all $x\in CI,f\in H_\ma$, there is a map $\tau_{x,f}$ such that
$\Phi(\lambda\theta_{x,f})=\tau_{x,f}(\lambda)\Phi(\theta_{x,f})$.
\end{lem}
\begin{proof}
For all $x\in H_\ma, f\in CI$, $\Phi(\theta_{x,f})=\theta_{y,g}$,
$f,g\in CI$. Therefore
$\Phi(\lambda\theta_{x,f})=\Phi(\theta_{\lambda
x,f})=\theta_{r(\lambda x),g}=\theta_{y,r'(\lambda f)}$ with $g\in
CI$. For a $e\in \varepsilon$, $\langle e,g\rangle r(\lambda
x)=\langle e, r'(\lambda f)\rangle y$. It follows from $g\in CI$
that $r(\lambda x)=\alpha y$, with $\alpha=\langle e,g\rangle
^{-1}\langle e,r'(\lambda f)\rangle$. Denote by
$\tau_{x,f}(\lambda)=\alpha$. Then
$\Phi(\lambda\theta_{x,f})=\tau_{x,f}\Phi(\theta_{x,f})$.
\end{proof}

\section{Main results}

\begin{thm}
 $\Phi:\mathcal{F}(H_\ma)\rightarrow \mathcal{F}(H_\ma)$ is a sujective rank one
 preserving additive map. $Im\Phi$ is neither contained in any $L_y$
 nor contained in any $R_g$ Then one of the following is true:
 (1) For all $x,f\in H_\ma$, $\Phi(\theta_{x,f})=\theta_{Ax,Cf}$
 where $A,C$ are injective quasi-modular maps on $H_\ma$;

 (2)For all $x,f\in H_\ma$, $\Phi(\theta_{x,f})=\theta_{Cf,Ax}$
 where $A,c$ are injective conjugate quasi-modular maps on $H_\ma$.

 \end{thm}
 \begin{proof}
 For any $x\in H_\ma$, there exists a $y\in H_\ma$ such that $\Phi(L_x)\subseteq
 L_y$. For all $f\in H_\ma$ with
 $\Phi(\theta_{x,f})=\theta_{y,C_xf}$.

 When $x$ in CI, and $y$ can be in CI, we have following claims.

 Claim 1. $C_xf$ is a map. In fact, putting $f_1=f_2$, we have
 $\Phi(\theta_{x,f_1}=\theta_{y,f_1}$,
 $\Phi(\theta_{x,f_2})=\theta_{y,C_xf_2}$. Since $y\in CI$, we get
 $C_x f_1=C_x f_2$.

Claim 2. $C_x$ is injective. Otherwise, there exists a $f_0\neq 0$
but $C_x f_0=0$. So $\Phi(\theta_{x,f_0})=\theta_{y,C_xf_0}=0$
contradicting to $\Phi$ preserving rank one.

Claim 3. $C_x$ is additive.
$\Phi(\theta_{x,f_1+f_2})=\theta_{y,C_x(f_1+f_2)}=\phi(\theta_{x,f_1})+\Phi(\theta_{x,f_2})=
\theta_{y,C_xf_1}+\theta_{y,C_xf_2}=\theta_{y,C_xf_1+C_xf_2}$. It
follows that $C_x(f_1+f_2)=C_xf_1+C_xf_2$ from $y\in CI$.

Claim 4. $C_x$ is a locally quasi-modular map. From Lemma 17,
$\Phi(\lambda\theta_{x,f})=\tau_{x,f}(\lambda)\Phi(\theta_{x,f})=\tau_{x,f}(\lambda)\theta_{y,C_x(f)}$,
$\Phi(\lambda\theta_{x,f})=\Phi(\theta_{x,\lambda^\ast
f})=\theta_{y,C_x(\lambda^\ast f)}$. So $C_x(\lambda^\ast
f)=\tau_{x,f}(\lambda)^\ast C_xf$.

Claim 5. $\tau_{x,f}$ is independent of $f$. In fact there are
$f,g\in CI$ with $C_x\bot C_xg$. So $C_x(\lambda f+\lambda
g)=\tau_{x,f+g}(\lambda)(C_xf+C_xg)=\tau_{x,f}(\lambda)C_xf+\tau_{x,g}(\lambda)C_xg$,
i.e.
$[\tau_{x,f+g}(\lambda)-\tau_{x,f}(\lambda)]C_xf-[\tau_{x,f+g}(\lambda)-\tau_{x,g}(\lambda)]C_xg=0$.
It follows that $[\tau_{x,f+g}(\lambda)-\tau_{x,f}(\lambda)]C_xf=0$
and $[\tau_{x,f+g}(\lambda)-\tau_{x,g}(\lambda)]C_xg=0$. Therefore
$\tau_{x,f}=\tau_{x,g}=\tau_{x,f+g}$.

For all $h\in CI$, $C_xh\bot C_xf$, $C_xh\bot C_xg$, we have
$\tau_{x,h}=\tau_{x,f+g}$. $\tau_{x,h}$ is multiplicative. In fact,
$\Phi(\lambda_1\lambda_2\theta_{x,h})=\tau_{x,h}(\lambda_1)\Phi(\lambda_2\theta_{x,g})
=\tau_{x,h}(\lambda_1)\tau_{x,h}(\lambda_2)\Phi(\theta_{x,h})$. And
since $\Phi(\theta_{x,h})=\theta_{y,C_x(h)}$, $y,C_x(h)\in CI$, we
have
$\tau_{x,h}(\lambda_1\lambda_2)=\tau_{x,h}(\lambda_1)\tau_{x,h}(\lambda_2)$.

For all $f=\sum\limits_i \alpha_i f_i$, $f_i\in CI$,
\begin{eqnarray*}
&&\Phi(\lambda \theta_{x,f})=\Phi(\lambda\theta_{x,\sum\limits_i
\alpha_if_i})=\Phi(\lambda\sum\limits_i \alpha_i^\ast
\theta_{x,f_i})\\
&=&\sum\limits_i\Phi(\lambda\alpha_i^\ast
\theta_{x,f_i})=\sum\limits_i\tau_{x,f+g}(\lambda\alpha_i^\ast)\Phi(\theta_{x,f_i})\\
&=&\sum\limits_i
\tau_{x,f+g}(\lambda)\tau_{x,f+g}(\alpha^\ast)\Phi(\theta_{x,f_i})\\
&=&\tau_{x,f+g}(\lambda)\Phi(\theta_{x,f})
\end{eqnarray*}
We infer $\tau_{x,f}(\lambda)=\tau_x(\lambda)$.

Claim 6. $\tau_x$ is a injective homomorphism from $\ma$ to $\ma$.
For all $\lambda_1,\lambda_2\in \ma$,
$\Phi((\lambda_1+\lambda_2)\theta_{x,f})=\tau_x(\lambda_1+\lambda_2)\theta_{y,C_xf}=
(\tau_x(\lambda_1)+\tau_x(\lambda_2))=\theta_{y,C_xf}$ and therefore
$C_x[\tau_x(\lambda_1+\lambda_2)^\ast
f-(\tau_x(\lambda_1)+\tau_x(\lambda_2))^\ast f]=0$. Since $C_x$ is
injective, we get $\tau_x(\lambda_1+\lambda_2)^\ast
f=[\tau_x(\lambda_1)^\ast+ \tau_x(\lambda_2)^\ast]f$. When putting
$f\in CI$ we have
$\tau_x(\lambda_1+\lambda_2)=\tau_x(\lambda_1)+\tau_x(\lambda_2)$.

From the above Claim we know $\tau_x$ is multiplicative.

If $\tau_x$ is not injective, then there exists a $0\neq\lambda_0\in
\ma$ but $\tau_x(\lambda_0)=0$. So
$\Phi(\lambda_0\theta_{x,f})=\tau_{x,f}(\lambda_0)\theta_{y,C_xf}=0$
which contradicting to $\Phi$ rank one preserving.

Claim 7. $\tau_x$ is independent of $x$.

Let $x_1,x_2\in CI$ with $\Phi(L_{x_1}^{CI})\subseteq L_{y_1}^{CI}$,
$\Phi(L_{x_2}^{CI})\subseteq L_{y_2}^{CI}$ and $\langle
y_1,y_2\rangle$, $y_1,y_2\in CI$.

For arbitrary $\lambda\in\ma$, we have
\begin{equation*}
\Phi(\lambda\theta_{x_1,f})=\tau_{x_1}(\lambda)\Phi(\theta_{x_1,f})=\tau_{x_1}(\lambda)\theta_{y_1,C_{x_1}f}
\end{equation*}
\begin{equation*}
\Phi(\lambda\theta_{x_2,f})=\tau_{x_2}(\lambda)\Phi(\theta_{x_2,f})=\tau_{x_2}(\lambda)\theta_{y_2,C_{x_2}f}
\end{equation*}
and so
\begin{equation}
\Phi(\lambda
\theta_{x_1+x_2,f})=\tau_{x_1}(\lambda)\theta_{y_1,C_{x_1}f}+\tau_{x_2}(\lambda)\theta_{y_2,C_{x_2}f}
=\tau_{x_1+x_2}(\lambda)\theta_{y_3,C_{x_1+x_2}f}
\end{equation}
When $\lambda=1$,
\begin{equation}
\theta_{y_1,C_{x_1}f}+\theta_{y_2,C_{x_2}f}=\theta_{y_3,C_{x_1+x_2}f}
\end{equation}
Since $\langle y_1,y_2\rangle=0,y_1,y_2\in CI$, $y_1\neq \alpha
y,y_2\neq \beta y$, for all $\alpha,\beta\in\ma,y\in H_\ma$. by
Corollary 6, there exists a invertible $\sigma\in\ma$ such that
$C_{x_1}f=\sigma C_{x_2}f$.

Similarly, it follows $C_{x_1}g=\nu C_{x_2}g$ from
$\Phi(\lambda\theta_{x_1,g})=\tau_{x_1}(\lambda)\theta_{y_1,C_{x_1}g}$,
$\Phi(\lambda\theta_{x_2,g})=\tau_{x_2}(\lambda)\theta_{y_2,C_{x_2}g}$
and putting $\lambda=1$.

Now consider the rank one map $\theta_{x_1+x_2,\lambda f+\lambda g}$
mapped by $\Phi$:
\begin{eqnarray*}
&&\Phi(\theta_{x_1+x_2,\lambda f+\lambda g})\\
&&=\tau_{x_1}(\lambda)\theta_{y_1,C_{x_1},f}+\tau_{x_1}(\lambda)\theta_{y_1,C_{x_1}g}+
\tau_{x_2}(\lambda)\theta_{y_2,C_{x_2}f}+\tau_{x_2}(\lambda)\theta_{y_2,C_{x_2}g}\\
&&=\tau_{x_1}(\lambda)\sigma\theta_{y_1,C_{x_2}f}+\tau_{x_1}(\lambda)\nu\theta_{y_1,C_{x_2}g}
+\tau_{x_2}(\lambda)\theta_{y_2,C_{x_2}f}+\tau_{x_2}(\lambda)\theta_{y_2,C_{X_2}g}\\
&&=\theta_{\tau_{x_1}(\lambda)\sigma
y_1+\tau_{x_2}(\lambda)y_2,C_{x_2}f}+\theta_{\tau_{x_1}(\lambda)\nu
y_1+\tau_{x_2}(\lambda)y_2,C_{x_2}g}
\end{eqnarray*}
On the other hand, by $\Phi$ preserving rank one, we know there
exist $y_3,h(\lambda)\in H_\ma$  such that
$\Phi(\theta_{x_1+x_2,\lambda f+\lambda
g})=\theta_{y_3,h(\lambda)}$. From Corollary again, there exist
$\alpha(\lambda),\beta(\lambda)\in\ma$ such that
\begin{equation}
\begin{cases}
&\tau_{x_1}(\lambda)\sigma
y_1+\tau_{x_2}(\lambda)y_2=\alpha(\lambda)y_3\\
&\tau_{x_1}(\lambda)\nu y_1+\tau_{x_2}(\lambda)y_2=\beta(\lambda)y_3
\end{cases}
\end{equation}
So
\begin{equation}
\begin{cases}
&\beta\tau_{x_1}(\lambda)\sigma
y_1+\beta\tau_{x_2}(\lambda)y_2=\alpha(\lambda)\beta(\lambda) y_3\\
&\alpha(\lambda)\tau_{x_x}(\lambda)\nu
y_1+\alpha(\lambda)\tau_{x_2}(\lambda)y_2=\alpha(\lambda)\beta(\lambda)
y_3
\end{cases}
\end{equation}
We infer
\begin{equation}
(\beta(\lambda)\tau_{x_1}(\lambda)\sigma-\alpha(\lambda)\tau_{x_1}(\lambda)\nu)y_1+(\beta(\lambda)\tau_{x_2}(\lambda)-
\alpha\tau_{x_2}(\lambda))y_2=0
\end{equation}
Since $\langle y_1,y_2\rangle=0$ and $y_1,y_2\in CI$, we get
\begin{equation}
\begin{cases}
&\beta(\lambda)\tau_{x_1}(\lambda)\sigma=\alpha(\lambda)\tau_{x_1}(\lambda)\nu\\
&\beta(\lambda)\tau_{x_2}(\lambda)-\alpha(\lambda)\tau_{x_2}(\lambda)=0
\end{cases}
\end{equation}
and we get $\beta(1)=\alpha(1)$, $\sigma=\nu$.

Thus
\begin{equation*}
C_{x_1}(\lambda
f)=\tau_{x_1}(\lambda)C_{x_1}f=\tau_{x_1}(\lambda)\sigma C_{x_2}f
\end{equation*}
\begin{equation*}
C_{x_1}(\lambda f)=\sigma C_{x_2}(\lambda f)=\sigma
\tau_{x_2}(\lambda)C_{x_2} f
\end{equation*}
and so
\begin{equation*}
\tau_{x_1}(\lambda)\sigma C_{x_2}f=\sigma\tau_{x_2}(\lambda)C_{x_2}f
\end{equation*}
which implying
\begin{equation*}
\tau_{x_1}(\lambda)=\tau_{x_2}(\lambda):=\tau(\lambda).
\end{equation*}

For every $x\in CI$, $\Phi(L_x^{CI})\subseteq L_y^{CI}$, for some
$y\in CI$. Then $y=\sum\limits_i \beta_i y_i$, where
$\{y_i\}\subseteq CI$ is a orthogonal set.
$\Phi(\theta_{x,f})=\theta_{y,C_xf}=\sum\limits_i \beta_i
\theta_{y_i,C_xf}$.

Since $\Phi$ is surjective, there are $\{x_i\}$ with
$\Phi(L_{x_i}^{CI})\subseteq L_{y_i}^{CI}$. For all $0\neq f\in CI$,
$\Phi(\theta_{x_i,f}+\theta_{x_0,f})=\theta_{y_i,C_{x_i},f}+\theta_{y_0,C_{x_0}f}$
is rank one, so that $C_{x_i}f=\gamma_i C_{x_0}f$ for some
invertible $\gamma_i\in\ma$. There are $\alpha_i\in\ma$ with
$\tau_{x_i}(\alpha_i)=\beta_i$,
$\tau_{x_i}(\delta_i)=\gamma_i^{-1}$. We consider
\begin{eqnarray*}
&&\Phi(\theta_{\sum\limits_i\delta_i\alpha_ix_i,f})\\
&=&\sum\limits_i\Phi(\theta_{\delta_i\alpha_ix_i,f})\\
&=&\sum\limits_i\tau_{x_i}(\delta_i\alpha_i)\Phi(\theta_{x_i},f)\\
&=&\sum\limits_i \gamma_i^{-1}\beta_i\theta_{y_i,C_{x_i}f}\\
&=&\sum\limits_i\gamma_i^{-1}\beta_i\theta_{y_i,\gamma_i C_xf}\\
&=&\sum\limits_i\beta_i\theta_{y_i,C_xf}
\end{eqnarray*}

It follows that
$\Phi(\theta_{x,f})=\Phi(\theta_{\sum\limits_i\delta_i\alpha_ix_i,f})$.
Since $\Phi$ preserving rank one and $0\neq f\in CI$, we get
$x=\sum\limits_i\delta_i\alpha_ix_i$.

For all $\lambda\in\ma$,
\begin{eqnarray*}
&&\Phi(\lambda\theta_{x,f})=\Phi(\lambda
\theta_{\sum\limits_i\delta_i\alpha_ix_i,f})\\
&=&\sum\limits_i\Phi(\lambda\delta_i\alpha_i\theta_{x_i,f})\\
&=&\sum\limits_i\tau_{x_i}(\lambda\delta_i\alpha_i)\Phi(\theta_{x_i,f})\\
&=&\sum\limits_i\tau(\lambda)\tau(\delta_i\alpha_i)\Phi(\theta_{x_i,f})\\
&=&\sum\limits_i\tau(\lambda)\Phi(\theta_{\delta_i\alpha_i x_i,f})\\
&=&\sum\limits_i\tau(\lambda)\Phi(\theta_{\sum\limits_i\delta_i\alpha_ix_i,f})\\
&=&\tau(\lambda)\Phi(\theta_{x,f})
\end{eqnarray*}
Thus we infer that $\tau_x$ is independent of $x$.

Claim 8. $C_x$ is independent of $x\in CI$. For $x_1,x_2\in CI$ with
$\Phi(L_{x_1}^{CI})\subseteq L_{y_1}^{CI}$,
$\Phi(L_{x_2}^{CI})\subseteq L_{y_2}^{CI}$, where $y_1,y_2\in CI$,
$y_1\neq\alpha y$, $y_2\neq\beta y$ for all $\alpha,\beta\in\ma,y\in
H_\ma$. Then we have
\begin{equation*}
\Phi(\theta_{x_1+x_2,f})=\theta_{y_1,C_{x_1}f}+\theta_{y_2,C_{x_2}f}
\end{equation*}
This yields $C_{x_1}f=\alpha_fC_{x_2}f$, for all $f\in H_\ma$. From
Lemma 9, we know there exists a $\alpha,\beta\in \ma$ such that
$C_{x_1}=\alpha C_{x_2}$, $C_{x_2}=\beta C_{x_1}$.

For arbatrary $x\in CI$ with $\Phi(L_x^{CI})\subseteq L_y^{CI},y\in
CI$ assume there exist $\alpha_0,\beta_0,\alpha_1,\beta_1\in\ma$,
$y_0,y_0'\in H_\ma$ such that $y=\alpha_0
y_0,y_1=\beta_0y_0;y=\alpha_1y_0',y_2=\beta_1y_0'$. We infer that
$\alpha_0,\alpha_1,\beta_0,\beta_1\in \ma$ are invertible from
$y,y_1,y_2\in CI$. Then we get $y_1=\gamma_1 y_0,y_2=\gamma_2 y_0$
for some $\gamma_1,\gamma_2\in \ma$ which contradicting to the
properties of $y_1,y_2$. Thus we have shown that either $y\neq
\alpha y',y_1\neq \beta y'$ or $y\neq \alpha y', y_2\neq\beta y'$
for all $\alpha,\beta\in\ma, y'\in H_\ma$. Consequently, $C_x$ and
$C_{x_1}$ differ only by a multiplicative $\alpha_x\in\ma$. By
absorbing this $\alpha_x$ in the first term of $\theta$, $C_x$
becomes independent of $x$. Denoting  $\alpha_xy$ by $A'x$, we get
\begin{equation}
\Phi(\theta_{x,f})=\theta_{y,C_xf}=\theta_{\alpha_xy,Cf}=\theta_{A'x,Cf}.
\end{equation}

Now for all $x\in H_\ma$ (may not in $CI$), $x=\sum\limits_i
\alpha_ix_i$, where $x_i\in CI$, with
$\Phi(\theta_{x_i,f})=\theta_{Ax_i,Cf}$. Then we have
\begin{equation}
\Phi(\theta_{x,f})=\Phi(\theta_{\sum\limits_i\alpha_ix_i,f})=\sum\limits_i\tau(\alpha_i)
\Phi(\theta_{x_i,f})=\sum\limits_i\tau(\alpha_i)\theta_{A'x_i,CF}=\theta_{\sum\limits_i\tau(\alpha_i)A'x_i,Cf}
\end{equation}
We denote $\sum\limits_i\tau(\alpha_i)Ax_i$ by $Ax$. Then for all
$x,f\in H_\ma$ we always have
\begin{equation}
\Phi(\theta_{x,f})=\theta_{Ax,Cf}
\end{equation}
When $f\in CI$, it is easy to show $A$ is a injective quasi-modular
map. And since $Ax $ is independent of $f$ we know that $A$ is
always a injective quasi-modular map.

Statement (2) can be shown similarly, and we can get our desired
results.

 \end{proof}

\bibliographystyle{amsplain}

\end{document}